\documentclass[final]{article}
\usepackage{showkeys}
\usepackage{amsmath, amssymb, epsf, a4wide,amssymb,amsmath,amsthm}
\usepackage{euscript}
\usepackage{graphics,graphicx}

\mathchardef\flat="115B

\newtheorem{thm}{Theorem}[section]

\newtheorem{cor}[thm]{Corollary}
\newtheorem{prop}[thm]{Proposition}

\begin{document}

\vskip 2.7cm

\centerline{\large{\bf A SEPARABLE MANIFOLD FAILING TO HAVE}}
\smallskip
\centerline{\large{\bf THE HOMOTOPY TYPE OF A CW-COMPLEX}}

\bigskip

%\centerline{ \small BY}

\centerline{ \normalsize ALEXANDRE GABARD}

\vskip 0.5cm

\smallskip
\newbox\abstract
\setbox\abstract\vtop{\hsize 11 cm \noindent

\footnotesize \noindent\textsc{Abstract.} We show that the
%classic
Pr\"ufer surface, which is a separable non-metrizable 2-manifold,
has not the homotopy type of a CW-complex. This will follow easily
from J. H. C. Whitehead's result: if
%you have
one has a good approximation of an arbitrary space by a
CW-complex, which fails to be a homotopy equivalence, then the
given space is not homotopy equivalent to
%any
a CW-complex.
%(Roughly if one fails, then everybody will fail.)
}

\centerline{\hbox{\copy\abstract}}

\bigskip

2000 {\it Mathematics Subject Classification.} {\rm 57N05, 57N99,
54D65, 57Q05.}

{\it Key words and phrases.} {\rm
%Manifolds,
Pr\"ufer manifold, Non-metrizable manifolds, Homotopy,
CW-complex.}

\normalsize

%\hskip 1.5cm

\section{Introduction}\label{sec1}

Our aim is to prove the following:

\begin{thm}\label{thm1}
The Pr\"ufer surface\footnote{Actually, the surface we consider is
not exactly the original Pr\"ufer surface, but rather
Calabi-Rosenlicht's slight modification of it.
%Pr\"ufer's example.
%the original Pr\"ufer
%surface.
%This will be clarified later in the paper.
} (which is an example of separable\footnote{A space is {\it
separable} if it has a countable dense subset.} non-metrizable
manifold\footnote{Here this means a Hausdorff topological space
which is locally Euclidean.}) has not the homotopy type of a
CW-complex.
\end{thm}

This might sound like a
%provocation, since Corollary 1 of Milnor's
%paper \cite{Milnor}, says that {\it every separable manifold has
%the homotopy type of a (countable) CW-complex.}
dissonance in view of Milnor's Corollary 1 \cite{Milnor}, which
says that {\it every separable manifold has the homotopy type of a
(countable) CW-complex.}

However the proof of Corollary 1 uses metrizability in a crucial
way. Remember that it works
%roughly
as follows: by results of J. H. C. Whitehead \cite{Whitehead}
%and \cite{Whitehead2}
it is enough to prove that our space is
dominated by a CW-complex. The first step is Hanner's theorem
\cite{Hanner}: a space that is locally an ANR is an ANR. Then
following Kuratowski, the space is embedded via $x \mapsto d(x,
\cdot)$ making
%``massive''
use of a (bounded) metric $d$ into a
Banach space
%(here already metrizability is needed)
as a closed subset of its convex hull $C$ (Wojdyslawski
\cite{Wojdyslawski}). Since it is an ANR, there is an open
neighborhood $U$ in $C$ retracting to our space. By transitivity
of domination, it is enough to prove that $U$ is dominated $f:P\to
U$ by a polyhedron $P$, which is constructed as the nerve of a
suitable cover. The ``submission'' map $g:U\to P$ is constructed
as the barycentric map attached to a partition of unity
(paracompactness is needed, but follows from metrizability).
%\cite{Stone},\cite{Rudin}).
Lastly the homotopy $fg\simeq 1_{U}$
comes from local convexity considerations. For a detailed
exposition see \cite{Palais}. (All this, being an elaboration of
the basic idea: embed
%your
the given space in an Euclidean space and triangulate an open
tubular neighborhood of it.)

%So one has to be a little bit careful with Milnor's Corollary 1 by
%applying it only to metrizable manifolds.
From this context it is quite clear that all manifolds in Milnor's
paper are implicitly assumed to be metrizable. In particular
Corollary 1 does not apply to the manifold constructed (under
CH=continuum hypothesis) by Rudin-Zenor \cite{RudinZenor}, which
is an example of hereditarily separable\footnote{Each subspace is
separable.} non-metrizable manifold. The question of the
contractibility of the Rudin-Zenor manifold then appears as an
interesting problem.

\section{The idea of the proof}\label{sec2}

Let us first give a loose description of the Pr\"ufer surface $P$.
%It can be thought of
We may think of $P$ as the (Euclidean) plane from which an
(horizontal) line has been suppressed, and then for each point of
the line a small bridge is introduced in order to connect the
upper half-plane $\cal H$ to the lower half-plane ${\cal
H}^{\sigma}$. (The formal construction of $P$ will be recalled in
$\S$ \ref{sec3}.)
%The only thing we have to know is that $P$ is
%separable.
For our argument, the only information which will be needed on $P$
is its separability. The idea
%is then simply that there is a
is to consider the natural map $f$ from the graph $K$ consisting
of two vertices and a continuum ${\frak c}=$(cardinality of ${\Bbb
R }$) of edges linking them, to $P$ given by going from a point in
$\cal H$ to a point in ${\cal H}^{\sigma}$ crossing through the
continuum of bridges at our disposal (see Figure 1).

\begin{figure}[h]
\centering
    \includegraphics[angle=0]{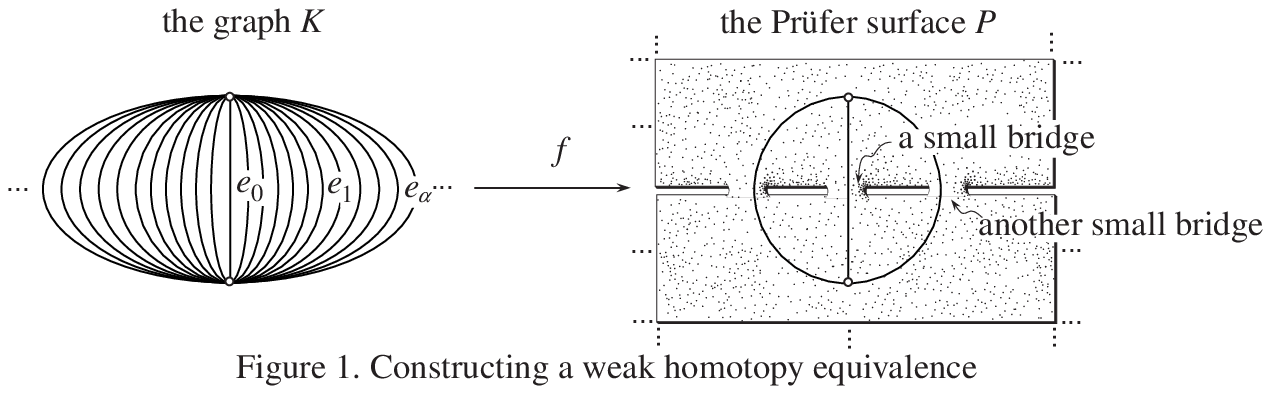}
\end{figure}

The fundamental group of $P$ is easily identified, via van
Kampen's theorem, as a free group on a continuum $\frak c$ of
%an uncountable set of
generators. Further, it is not hard to check that the higher
homotopy groups of $P$ vanish, i.e. $\pi_i(P)=0$ for all $i\ge 2$.
(Details will be given in $\S$ \ref{sec4}.) It follows that the
map $f:K \to P$ is a weak homotopy equivalence\footnote{i.e.
induces isomorphisms on all homotopy groups $\pi_i$.}, which turns
out {\it not} to be (as we will soon explain) a strong homotopy
equivalence. Then Whitehead's Theorem 1 \cite{Whitehead} implies
that $P$ has not the homotopy type of a CW-complex. To explain why
$f$ is not a homotopy equivalence, we first observe that if there
were a homotopy inverse $g:P \to K$, then
%since $g$ induces an
%isomorphism on the fundamental groups $\pi_1$, the map
$g$ has to be onto. [This, because if $g$ misses a point $p\in K$,
then $g$ factors through $K-p$, and so by functoriality the
inclusion $K-p \hookrightarrow K$ has to induce an epimorphism on
the $H_1$ (first homology group). But our graph $K$ is easily
verified to be such that
%that the suppression of
for any point $p\in K$
%leads to an
the inclusion $K-p \hookrightarrow K$
%, which
fails to induce a epimorphism on the $H_1$]. But then pulling-back
the uncountable collection of open 1-cells $(e_\alpha)_{\alpha \in
{\Bbb R}}$ of $K$, we get in $P$ an uncountable family
$(g^{-1}(e_\alpha))_{\alpha \in {\Bbb R}}$ of pairwise disjoint
open sets. This is a contradiction, since $P$ is separable. \qed

\section{Construction of the Pr\"ufer surface $P$}\label{sec3}
\medskip

The following construction is due to Pr\"ufer, first described in
print by Rad\'o \cite{Rado}. (As it is well-known, this was also
the ``first'' example of a non-triangulable surface,
%and thereby
%which
and played an important role in clarifying foundational aspects of
Riemann surfaces theory.)
% Now, the formal definition.
We use ${\Bbb C}$ as model for the Euclidean plane. The idea is to
consider the set $P_{0}$ formed by the (open) upper half-plane
${\cal H} = \{ z : {\rm Im }(z)>0 \}$ together with the set of all
rays emanating from point of ${\Bbb R}$ and pointing out in the
upper half-plane. Then we topologize $P_0$ with the usual topology
for $\cal H$, and by taking as neighborhoods of a point $r$ which
is a ray (say emanating from $x\in {\Bbb R}$) an (open) sector of
rays deviating by at most $\varepsilon$ radian from $r$, together
with the points of $\cal H$ between the two rays and at
(Euclidean) distance smaller than $\varepsilon$ from $x$ (see
Figure 2).

\begin{figure}[h]
\centering
    \includegraphics[angle=0]{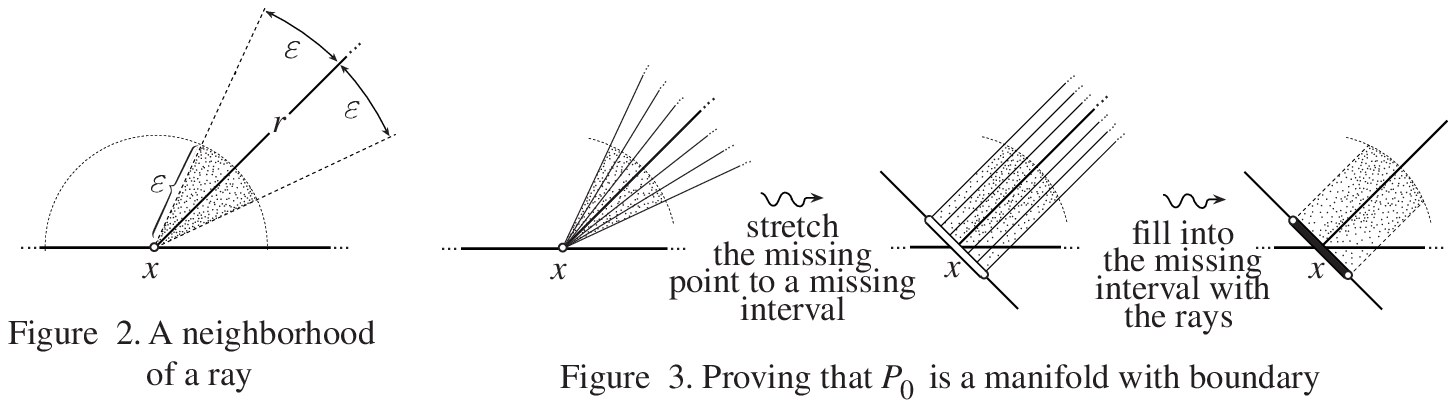}
\end{figure}

The space $P_0$ is a surface-with-boundary, as it is easily seen
by following the pictures from Figure 3. Observe that $P_0$ has a
continuum ${\frak c}$ of boundary components each homeomorphic to
the real line ${\Bbb R}$.

Now given a manifold-with-boundary $W$, there are
%(at least)
two obvious ways to obtain a manifold $M$ (without boundary): a
first method is by {\it collaring}, set $M=W \cup_{{\rm
id_{\partial W}}} (\partial W \times [0,1))$ (glue $W$ with the
cylinder having the boundary $\partial W$ as basis, along their
boundaries), and a second option is by {\it doubling}, $M=W
\cup_{{\rm id_{\partial W}}} W$ (two copies of $W$ are glued along
their boundaries).

For $W=P_0$, the process of collaring leads to the ``original''
Pr\"ufer surface (the one described in \cite{Rado}). In this case
%an uncountable family of half-plane is glued, so
there is in $M$ an uncountable family of pairwise disjoint open
sets. This implies that $M$ is non-separable. (And so it did not
really interest us: remember our purpose is to point out that the
metrizability hypothesis in Milnor's Corollary 1 is essential.)

The second option leads to the surface $P$ we are interested in
since it is separable. We call it also the Pr\"ufer surface (even
though
%to my knowledge
it seems to appear explicitly only in the paper by
Calabi-Rosenlicht \cite{Calabi}).

%Finally the Pr\"ufer surface $P=P_0 \cup_{{\rm id}} P_0$ is
%defined by taking two copies of $P_0$, and gluing them along their
%boundaries.

\begin{prop}
  The Pr\"ufer surface $P$ obtained by the process of doubling, is a connected (Hausdorff) $2$-manifold which is separable,
  but contains an uncountable discrete subspace
  (and therefore is non-metrizable).
\end{prop}

\proof Observe that the rational points ${\Bbb Q}+i {\Bbb Q_{>0}}$
give a countable dense subset of $P_0$, and so $P$ is clearly
separable. Further notice that the set of all rays $(r_x)_{x\in
{\Bbb R}}$, say orthogonal to ${\Bbb R}$ gives an uncountable
discrete subset of $P$, since an open neighborhood of such a ray
cuts out only this single one from the whole family.
%Finally this implies
It follows that $P$ is not hereditarily separable, and so not
second countable, and therefore non-metrizable. (As it is
well-known metrizability and second countability are equivalent
for (connected) manifolds. Actually in our situation since $P$ is
separable, the non-metrizability of $P$ can also be deduced from
the elementary fact that metrizable plus separable imply second
countable.)
\endproof

At this stage one could already observe the following:

\begin{cor}
  The Pr\"ufer surface $P$ (and more generally any non-metrizable manifold) is not homeomorphic to a CW-complex.
\end{cor}

\proof This follows from the fact proved by Miyazaki
\cite{Miyazaki} that CW-complexes are always paracompact, and the
equivalence between the concepts of paracompactness and
metrizability, when spaces are restricted to be manifolds.
%for manifolds.
\endproof

\section{ Homotopy groups
%computation
of the Pr\"ufer surface}\label{sec4}
\medskip

We now investigate the homotopy groups of $P$. First, the
fundamental group $\pi_1$:

\begin{prop}
  $\pi_1(P)$ is a free group on a continuum $\frak c$ of generators.
\end{prop}

\proof (Following \cite{Baillif}). For all $x\in \Bbb R$, let
$U_x$ be the open neighborhood depicted in Figure 2 with $r$
chosen orthogonal to ${\Bbb R}$ and $\varepsilon={\pi \over 2}$.
Let then $B_x$ be $U_x$ taken together with its symmetrical copy
$U_x^{\sigma}$, so $B_x=U_x\cup U_x^{\sigma}$ is an open set of
$P$ (we can think of it as a ``bridge'' linking the upper to the
lower half-plane).
%Consider now the following cover of
%$P$:
For all $x\in \Bbb R-\{ 0 \}=\Bbb R^{\ast}$, put ${O}_x={\cal H}
\cup {\cal H}^{\sigma} \cup B_0 \cup B_x$. The collection
$(O_x)_{x\in{\Bbb R}^{\ast}}$ form an open cover of $P$, which
satisfies van Kampen's theorem hypothesis, since $O_x \cap O_y$ is
arcwise-connected. Furthermore $\cap_{x\in{\Bbb R}^{\ast}}
O_x={\cal H} \cup {\cal H}^{\sigma} \cup B_0$, this being
homeomorphic to $\Bbb C - \Bbb R$ union an open interval from
$\Bbb R$ (by the same kind of argument as the one cartooned in
Figure 3), and so in particular simply connected. Moreover each
member $O_x$ of this cover is homeomorphic to $\Bbb C - \Bbb R$
union two disjoint (real) intervals, and so
%belong to the
has the homotopy type of the circle $S^1$. The result follows by
van Kampen's theorem.
\endproof

This proof gives not only the abstract algebraic structure of the
group $\pi_1(P)$, but also a concrete geometric description for
generators. They can be chosen as the loops going first down from
$i=\sqrt{-1}$ to $-i$ straightforwardly, and then going up by
travelling along an elliptic arc crossing the real axis through
some non-zero real number $x$ and turning back again to $i$ (see
Figure 1). (One further has to agree that all times we cross the
missing real axis at some abscissae $x\in {\Bbb R}$, we cross it
through the ray orthogonal to $\Bbb R$ emanating from $x$.)
Tautologically, it follows that the map $f:K \to P$, being
essentially defined in the same way, induces an isomorphism on the
fundamental groups.

Next the higher homotopy groups are simply given by:

\begin{prop}
  $\pi_i(P)=0$ for all $i\ge 2$.
\end{prop}

\proof We consider the open cover $\cal U$ of $P$ given by
%the following list of open subsets
$\cal H$, ${\cal H}^{\sigma}$ and for all $x \in \Bbb R$ the
bridge $B_x=U_x\cup U_x^{\sigma}$ (previously defined). Let now
$\gamma \in \pi_i(P)$, and choose $c: S^i \to P$ a representing
map. Pulling back the cover $\cal U$ by the map $c$, we see by
compactness of the sphere $S^i$ that there is a finite subset
$\cal F$ of $\cal U$ such that the image $c(S^i)$ is contained in
$\cup_{U \in {\cal F}} U =:V$. But it is clear (again by Figure~3)
that $V$ is homeomorphic to $\Bbb C - \Bbb R$ union a finite
number of open intervals (if not something even simpler like only
$\cal H$ or $\cal H$ union a finite number of bridges). In any
case the space $V$ has the homotopy type of a wedge of circles
(maybe an empty one), and so has vanishing $\pi_i$ for $i\ge 2$.
We are done.
\endproof

The proof of Theorem 1.1 is then complete. \qed

\section{The case of non-Hausdorff manifolds}\label{sec5}
\medskip
%Concerning the general question of knowing under which
%circumstances a manifold has the homotopy type of a CW-complex,
%notice that the Hausdorff condition included in the definition of
%a manifold is quite essential.

We conclude by making simple observations concerning complications
arising in the relation between manifolds and CW-complexes,
%becomes much more troublesome,
if we drop Hausdorff from the definition of a manifold.
%Without it,
Then already one of the simplest example of
%non-Hausdorff
``manifold'', the so called {\it line with two origins} (obtained
from two copies of ${\Bbb R}$ by identifying corresponding points
outside the origin, see Figure 4) fails to have the homotopy type
of a CW-complex (and this in spite of the fact that it is
well-behaved from the point of view of second countability).
Actually we even have a worse situation:

\begin{prop}
  The line with two origins $R$ has not the homotopy type of any Hausdorff topological space.
\end{prop}

\begin{figure}[h]
\centering
    \includegraphics[angle=0]{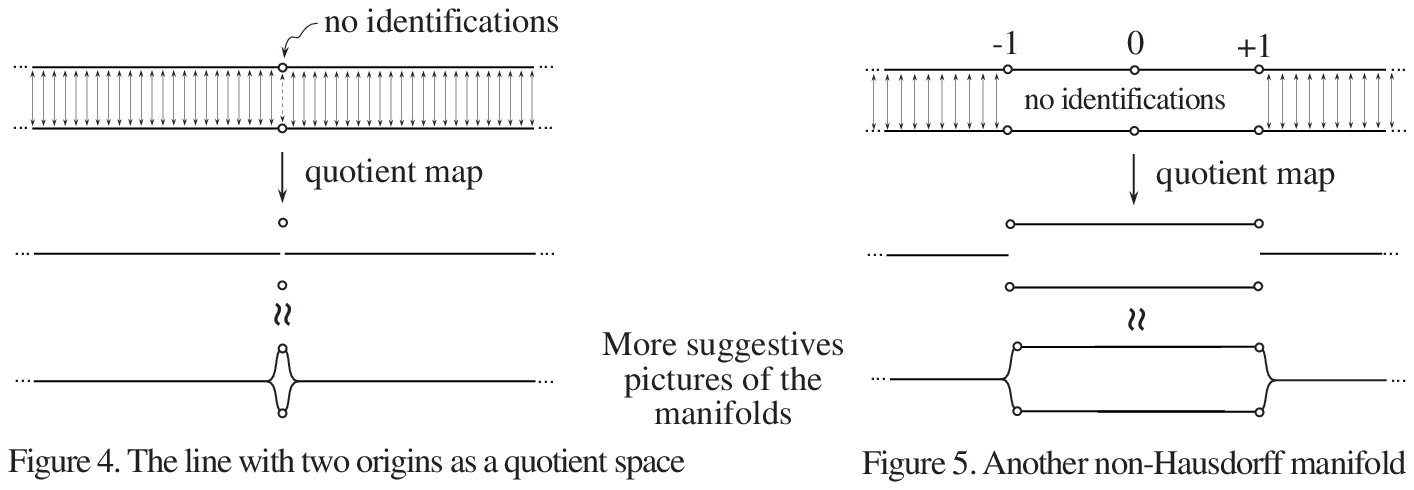}
\end{figure}

\proof We need
%to make
two preliminary remarks.

$\bullet$ First remember that there is a general Hausdorffization
process applicable to any space $X$, which leads to an Hausdorff
space $X_{\rm Haus}$ with a map $X \to X_{\rm Haus}$. (This is by
moding out the given space by the equivalence relation generated
by the inseparability relation.) It has the property that any
(continuous) map from $X$ to an Hausdorff space $H$ factors
through $X_{\rm Haus}$.

$\bullet$ Second by Mayer-Vietoris it is easy to check that
$H_1(R)\cong \Bbb Z$.

We are now in position to prove 5.1. Assume there is a homotopy
equivalence $f: R \to H$ between $R$ and some Hausdorff space $H$.
Then $f$ factors through $R_{\rm Haus}$, which is nothing else
than the usual real line $\Bbb R$. But this being contractible, it
follows by functoriality that the morphism $H_1(f)$ is zero, in
contradiction to the non-vanishing of $H_1(R)$. \endproof

Furthermore it is not difficult to compute the homotopy groups of
$R$ (for example by looking at $R$ as the leaf space of $\Bbb C
-0$ foliated by vertical lines, and then applying the exact
homotopy sequence of a fibration). We conclude that $R$ is an
Eilenberg-Mac Lane space of type $K({\Bbb Z},1)$ not homotopy
equivalent to the circle $S^1$.

Finally, a variant of the line with two origins (see Figure 5)
%obtained not by
%doubling only the origin, but a closed interval of some positive
%length,
leads to a non-Hausdorff manifold  which is easily seen to be
homotopy equivalent to the circle $S^1$. So,
%that
it is not non-Hausdorffness as a rule, that leads us outside the
%world
class $\cal W$ of spaces having the homotopy type of a CW-complex,
but much more
%much
%more the fact that non-Hausdorff manifolds
%may present geometric behavior of very narrow ``bifurcations''.
%more
the strange geometric behavior of ``extremely narrows
bifurcations'' presented by some non-Hausdorff manifolds, which
%is
appears as something
%too ``instantaneous'' and
alien to the combinatorial
nature of CW-complexes.
%to be restored
%build
%catched
%recorded
%displayed
%represented
%restored
%by a in essence combinatorial object (CW-complex).
%by a CW-complex

\bigskip \noindent{\bf Acknowledgments.}
%It is a pleasure to thank Mathiew Baillif for
%%transmitting perceptions on some
%precious guidance in the theory of non-metrizable manifolds.
I am
%owe
very much
%indebted
obliged to Matthew Baillif for precious guidance through some
theory of non-metrizable manifolds, and to
%wish also to thank
Andr\'e
Haefliger for suggesting the method
%giving the idea how
to compute $\pi_i(R)$.

\bigskip

\newbox\adress
\setbox\adress\vtop{\hsize 14.5 cm \noindent

\footnotesize
%%M. Baillif
%\noindent \textsc{Universit\'e de Gen\`eve, Section de
%Math\'ematiques, 2-4, rue du Li\`evre, CP 240, 1211 Gen\`eve 24,
%Suisse}

%\noindent{\it E-mail address:} \verb"baillif@math.unige.ch"

%A. Gabard
\medskip
\noindent \textsc{Universit\'e de Gen\`eve, Section de
Math\'ematiques, 2-4, rue du Li\`evre, CP 240, 1211 Gen\`eve 24,
Suisse}

\noindent{\it E-mail address:} \verb"alexandregabard@hotmail.com"}

\vskip -0.0 cm

\centerline{\hbox{\copy\adress}}

\end{document}